\documentclass[12pt]{article}

\usepackage[utf8]{inputenc}
\usepackage[T1]{fontenc}
\usepackage[english]{babel}
\usepackage{amsmath,amssymb,amsthm}

\title{Finite Decompositions for Initial Baire‑1 Subclasses and Their Characterization via Homogeneous Sequences}
\author{Alexey Ostrovsky\\
\texttt{detefu6@gmail.com}}

\date{}

\usepackage{amsthm}

\newtheorem{theorem}{Theorem}

\begin{document}
\maketitle

\begin{abstract}
We announce three results concerning initial subclasses of the Baire--1 functions.
First, we obtain a decomposition into two continuous maps for a basic initial subclass.
Second, we give a pointwise description of initial subclasses in terms of point-homogeneous sequences.
Third, we present a decomposition of maps into a finite number of open and closed maps.
\end{abstract}

For almost two centuries after the classical pointwise definition of discontinuity was introduced, no finer pointwise descriptions existed that could distinguish individual subclasses in the Baire classification of discontinuous functions. 

 Until now, the question of pointwise descriptions for subclasses of the first Baire class had not even been formulated. It naturally arises from a different problem considered below: which Baire‑1 functions admit a finite decomposition into continuous pieces? The classical theory of countable decompositions is of no use here, since complete metrizability is only finitely additive, and no classification existed that would distinguish functions by the minimal number of continuous components. The present work introduces the missing finite framework and begins a systematic study of pointwise structures suitable for establishing finite decomposability.

\section{Finite Decompositions: Main Ideas}

The problem of decomposing a function into a finite number of
continuous components was first considered in my earlier work
on finitely continuous functions \cite{OstrovskyFinitelyContinuous}.

The main result established there is that every $O \cap F$-measurable function
can be decomposed into three continuous functions.

As an immediate consequence, every one-to-one function $f$ such that both
$f$ and $f^{-1}$  are $O \cap F$-measurable can be decomposed into nine homeomorphisms.

The minimal possible number
of homeomorphisms in such a decomposition is not known.

In the
special case of $O \cap F$-measurable functions, the optimal number
turns out to be four. For other classes, nothing is known at present.

This result shows that the initial classes admit finite continuity bounds, and motivates the introduction of the higher classes \(C_n\).

\begin{theorem}
Every $O\cap F$-measurable map and every $O\cup F$-measurable map
$f:X\to Y$ is $2$-continuous.
\end{theorem}

To treat these remaining cases, we introduce the $C_n$-classification
of sets. Only the initial classes will be described here; the general scheme is given in Appendix A.

\section{The $C_n$-Classification of Sets (Initial Classes)}

The classes $C_n$ are defined by applying $n$ successive operations
of union or intersection to open sets $O$ and closed sets $F$.

\subsection*{Notation}

The first class $C_{1}$ consists of the two basic types of sets

\[
C_{1}=\{LC_{1},\, F\cup O\},
\]

where $LC_{1}=O\cap F$ denotes the sets obtained by one intersection of an
open and a closed set, and $F\cup O$ denotes the sets obtained by one union
of an open and a closed set.

In the sequel we write $LC_{n}$ for the union of $n$  sets of type $LC_{1}.$

For higher classes we use the mixed notation

\[
F \cup LC_n, \qquad O \cup LC_n,
\]

which denotes the sets obtained by one additional union of $F$ or $O$ with a
set from the class $LC_n$.

Thus, for example,

\[
C_2 = \{LC_1 \cup F,\; LC_1 \cup O\},
\]

and

\[
C_3 = \{LC_2,\; F \cup LC_1 \cup O\}.
\]

These notations will be used throughout the paper.

\medskip
We list the first few classes explicitly.

\textbf{Class $C_0$.} Open sets $O$ and closed sets $F$.

\textbf{Class $C_1$.} Sets of the form $LC_1$ and $F \cup O$.

\textbf{Class $C_2$.} Sets of the form $LC_1 \cup F$ and $LC_1 \cup O$.

\textbf{Class $C_3$.} Sets of the form $LC_2$ and $F \cup LC_1 \cup O$.

\textbf{Class $C_4$.} Sets of the form $LC_2 \cup F$ and $LC_2 \cup O$.

\textbf{Class $C_5$.} Sets of the form $LC_3$ and $F \cup LC_2 \cup O$.

\section{Decomposition into Open and Closed Maps}

In this section we consider continuous maps between metric spaces and
discuss several decomposition results for open--constructible and closed--constructible
maps.  Here a continuous map $f:X\to Y$ is called open–constructible (respectively,
closed–constructible) if the image of every open set  (respectively, every closed  set) under $f$ is an open–constructible
set (respectively, a closed–constructible set), that is, a set belonging to one of the
classes $C_{n}$ generated from open and closed sets by finitely many elementary
operations of union and intersection.

It was shown in \cite{OstrovskyClosedConstructible} that every continuous
closed--$LC_n$ map is countably closed; equivalently, such a map can be represented
as a countable union of closed maps.  
However, the analogous statement for open--constructible maps fails: as demonstrated
by \cite{OsipovRemarksOstrovsky}, there exists a continuous open--$C_1$
map which is not countably open.

Recently  it was established that every continuous
open--$LC_1$ map, as well as every continuous open--$(O \cap F)$ map, admits a
decomposition into one open map and one closed map \cite{OstrovskyNafpaktos2023}.  
We extend this observation to several broader classes.

\begin{theorem}
Every continuous open--$(LC_1 \cup F)$ map $f : X \to Y$ can be decomposed
into the union of one open map and two closed maps.
\end{theorem}

\begin{theorem}
Every continuous open--$(LC_1 \cup O)$ map $f : X \to Y$ can be decomposed
into the union of one closed map and two open maps.
\end{theorem}

\begin{theorem}
Every continuous open--$LC_2$ map $f : X \to Y$ can be decomposed into the
union of two closed maps and two open maps.
\end{theorem}

\medskip

\noindent\textbf{Question.}
If the image of every open set is a $C_n$ set, is $f$ the union of countably
many closed and open maps?

\section{The Sets $D_0,D_1,D_2,D_3$ and the Canonical Compacts $S_1,S_2,S_3$}

We use the following notation:

\[
D_0 = \{0\}, \quad
D_1 = \left\{ \frac{1}{2^{n_1}} : n_1 \in \mathbb{N}^+ \right\}, \quad
D_2 = \left\{ \frac{1}{2^{n_1}} + \frac{1}{2^{n_2}} : n_1,n_2 \in \mathbb{N}^+ \right\}.
\]

The sets $D_3$ are defined analogously by triples of indices and will be used
in the description of the canonical compact $S_3$.

Using these sets, we define the canonical compacts:

\[
S_1 = D_0 \cup D_1,\quad
S_2 = D_0 \cup D_1 \cup D_2,\quad
S_3 = D_0 \cup D_1 \cup D_2 \cup D_3.
\]

These canonical compacts will be used to describe pointwise measurability and
the discontinuity index of a map.

\section{Basic Properties of the Sets $D_0, D_1, D_2, D_3$}

We list several easily verifiable properties of the sets $D_0, D_1, D_2, D_3$
within the structures $S_1$, $S_2$, and $S_3$.

\medskip

\textbf{In $S_1$:}

\[
D_0 \text{ is an } F\text{-set but not an } O\text{-set in } S_1,
\]

\[
D_1 \text{ is an } O\text{-set but not an } F\text{-set in } S_1.
\]

\medskip

\textbf{In $S_2$:}

\[
D_0 \cup D_2 \text{ is an } F \cup O\text{-set but not an } LC_1\text{-set in } S_2,
\]

\[
D_1 \text{ is } LC_1 \text{ but not } F \cup O \text{ in } S_2.
\]

\medskip

\textbf{In $S_3$:}

\[
D_0 \cup D_2 \text{ is } F \cup LC_1 \text{ but not } O \cup LC_1 \text{ in } S_3,
\]

\[
D_1 \cup D_3 \text{ is } O \cup LC_1 \text{ but not } F \cup LC_1 \text{ in } S_3.
\]

\section{Characterization of Non-$C_2$-Measurable Maps}

In this section we describe the local structure of maps that are not
$C_2$-measurable at a point. .

\subsection*{Notation}

Let $S_n$ be the canonical compact defined in Section~4.
For a point $x\in X$, let $S_n(x)$ denote the image of $S_n$ under a
homeomorphism sending $0$ to $x$. Thus

\[
S_n(x)=D_0(x)\cup D_1(x)\cup\dots\cup D_n(x),
\]

where $D_j(x)$ is the image of the layer $D_j$ of $S_n$.
Each $D_j(x)$ consists of points of the form $x_{i_1\ldots i_j}$,
but we suppress indices since all branches are homeomorphic.

\medskip

We now characterize non-$C_2$-measurability at a point.

\subsection*{Non-$F\cup LC_1$-measurability at a point}

A map $f:X\to Y$ is said to be not $F\cup LC_1$-measurable if there exists
an open set $U\subset Y$ such that $f^{-1}(U)$ is not an $F\cup LC_1$ set.

A map $f$ is not $F\cup LC_1$-measurable at a point $x\in X$ if there exist
an open neighbourhood $U\subset Y$ of $y=f(x)$ and the canonical compact

\[
S_3(x)=D_0(x)\cup D_1(x)\cup D_2(x)\cup D_3(x),
\]

such that

\[
f^{-1}(U)\cap S_3(x)=D_1(x)\cup D_3(x)
   =\{x_i:i \in\mathbb N\}\cup\{x_{i,j,k}:i,j,k\in\mathbb N\}.
\]

A map is not $F\cup LC_1$-measurable if and only if it is not
$F\cup LC_1$-measurable at some point.

\subsection*{Non-$O\cup LC_1$-measurability at a point}

Similarly, a map $f$ is not $O\cup LC_1$-measurable at a point $x\in X$
if there exist an open neighbourhood $U\subset Y$ of $y=f(x)$ and the same
canonical compact $S_3(x)$ such that

\[
f^{-1}(U)\cap S_3(x)=D_0(x)\cup D_2(x)
   =\{x\}\cup\{x_{i,j}:i,j\in\mathbb N\}.
\]

A map is not $O\cup LC_1$-measurable if and only if it is not
$O\cup LC_1$-measurable at some point.

\section{Non-$C_3$-Measurability at a point}

\textbf{Definition.}
A map $f : X \to Y$ is not $LC_2$-measurable if there exists
an open set $U \subset Y$ such that $f^{-1}(U)$ is not an
$LC_2$ set.

We say that $f$ is not $LC_2$-measurable at a point $x \in X$
if there exists an open neighbourhood $U \subset Y$ of
$y = f(x)$ and a canonical compact  $S_4(x)$ such that

\[
f^{-1}(U) \cap S_4(x)
   = \{x\} \cup \{x_{i j}:  i, j \in \mathbb{N}\} \cup \{x_{i j k l}: i, j, k,l \in \mathbb{N}\}.
\]

The map $f$ is not $LC_2$-measurable if and only if it is not
$LC_2$-measurable at some point.

Analogously, $f$ is not $F \cup LC_1 \cup O$-measurable at a point $x$
if there exist an open set $U \subset Y$ and a 
$S_4(x)$ such that

\[
f^{-1}(U) \cap S_4(x)
   =\{x_{i} : i \in \mathbb{N}\}   \cup \{x_{i j k} : i, j, k \in \mathbb{N}\}.
\]

Similarly, $f$ is not $F \cup LC_1 \cup O$-measurable if and only if it is not
$F \cup LC_1\cup O$-measurable at some point.

\medskip


\medskip

\appendix

\section{General Pointwise Scheme for the Classes $C_n$}

\subsection*{General Definition of the Classes $C_k$}

For completeness we give the formal definition of the hierarchy $\{C_k\}$.

Let $O$ denote the family of open sets and $F$ the family of closed sets.
Define

\[
C_0 = \{O, F\}.
\]

For $k \ge 1$, the class $C_k$ consists of all sets obtained from $C_{k-1}$
by one additional elementary operation of union or intersection with $O$ or $F$:

\[
C_k = \{\, A\cup O,\; A\cup F,\; A\cap O,\; A\cap F : A\in C_{k-1}\,\}.
\]

Thus each level $C_k$ is generated by applying one elementary operation to
the sets of the previous level.

\medskip
\noindent

\noindent
The hierarchy $\{C_k\}$ is generated by applying two elementary operations
(open and closed) at each step. Only this recursive structure will be used
below; no further internal description is needed.

\medskip
\noindent
\textbf{Even and odd classes.}
The two pointwise obstructions at each level $C_k$ arise from selecting
either the even layers or the odd layers of $S_{k+1}(x)$:

\[
\text{even part: } D_0(x)\cup D_2(x)\cup\dots,\qquad
\text{odd part: } D_1(x)\cup D_3(x)\cup\dots.
\]

This is why every class $C_k$ has exactly two canonical pointwise forms:
one built from the even components $D_{2j}(x)$ and one built from the
odd components $D_{2j+1}(x)$.

The hierarchy $\{C_n\}$ is determined by pointwise non-$C_n$-measurability.
Each level corresponds to one of four finite pointwise decompositions of the
canonical compact $S_m(x)$ into its components $D_j(x)$.

\subsection*{Odd levels $C_{2k+1}$}

Non-$C_{2k+1}$-measurability at a point $x$ has two forms.

\smallskip
\textit{(1) Not $LC_{k+1}$-measurable at $x$.}
There exists an open neighbourhood $U$ of $f(x)$ such that

\[
f^{-1}(U)\cap S_{2k+2}(x)
   = D_0(x)\cup D_2(x)\cup\dots\cup D_{2k+1}(x).
\]

\textit{(2) Not $F\cup LC_k\cup O$-measurable at $x$.}
There exists an open neighbourhood $U$ of $f(x)$ such that

\[
f^{-1}(U)\cap S_{2k+2}(x)
   = D_1(x)\cup D_3(x)\cup\dots\cup D_{2k+1}(x).
\]

\subsection*{Even levels $C_{2k}$}

Non-$C_{2k}$-measurability at a point $x$ also has two forms.

\smallskip
\textit{(3) Not $LC_k\cup F$-measurable at $x$.}
There exists an open neighbourhood $U$ of $f(x)$ such that

\[
f^{-1}(U)\cap S_{2k+1}(x)
   = D_1(x)\cup D_3(x)\cup\dots\cup D_{2k+1}(x).
\]

\textit{(4) Not $LC_k\cup O$-measurable at $x$.}
There exists an open neighbourhood $U$ of $f(x)$ such that

\[
f^{-1}(U)\cap S_{2k+1}(x)
   = D_0(x)\cup D_2(x)\cup\dots\cup D_{2k}(x).
\]

\medskip

Thus every class $C_n$ is generated by a unique finite decomposition of the
compact $S_m(x)$ into its even or odd components.
These four patterns describe all possible local obstructions to
$C_n$-measurability.
\section{Point-Homogeneous Sequences}

Throughout the paper we used various canonical compacts

\[
S_n = D_0 \cup D_1 \cup \dots \cup D_n,
\]

as well as their even and odd substructures

\[
D_0 \cup D_2 \cup \dots \cup D_{2k},
\qquad
D_1 \cup D_3 \cup \dots \cup D_{2k+1},
\]

and similar finite unions of the sets $D_m$.

Although these objects appear in different contexts — in the description of
pointwise $C_n$-measurability, in the definition of discontinuity index, and in
the characterization of non-$C_2$-measurable maps — they are all instances of a
single general notion: \emph{point-homogeneous sequences}.

\medskip

Let \(I = [0,1]\). A sequence \(S \subset I\) is called point-homogeneous  if it has exactly one 
limit point \(p \in I\), and if the sets …

\[
U_n \cap (S \setminus \{p\})
\]

are pairwise homeomorphic for some countable base $\{U_n\}$ at $p$.

\medskip

In this sense, the canonical compacts $S_n$ are natural generalizations of an
ordinary convergent sequence with a single limit point.  They preserve the same
local symmetry at the limit point, but allow arbitrarily rich finite branching
structures.  This makes them ideally suited for describing pointwise behaviour
of maps in the hierarchy of the classes $C_n$.

Thus, all constructions used in the paper — even and odd unions of the sets
$D_m$, the compacts $S_n$, and their substructures — are unified by the single
concept of point-homogeneity.

\section{ Finite Discontinuity Index at a point}

Let $x$ be a point of discontinuity of a map $f$. If $f$ is $C_{1}$-measurable
at $x$, then we say that $x$ has finite discontinuity index $1$.

By induction, if $f$ is not $C_{n}$-measurable at $x$, but is
$C_{n+1}$-measurable at $x$, then $x$ is said to have finite discontinuity
index $n+1$.

\noindent
Thus, the non-$C_{n}$ behaviour of a map is captured by points at which the
discontinuity index is a finite integer determined by the associated
point-homogeneous sequences.


\medskip


\begin{thebibliography}{99}

\bibitem{OstrovskyFinitelyContinuous}
A.~Ostrovsky,
\textit{Finitely continuous functions},
Topology and its Applications,
vol.~261,
2019,
pp.~46--50,
doi:10.1016/j.topol.2019.05.004.


\bibitem{OstrovskyClosedConstructible}
A.~Ostrovsky,
\textit{Closed-constructible functions},
Topology and its Applications,
vol.~261,
2019,
pp.~46--50,
doi:10.1016/j.topol.2019.05.004.

\bibitem{OsipovRemarksOstrovsky}
A.~Osipov,
\textit{Remarks on Ostrovsky's theorem},
Sibirskie Elektronnye Matematicheskie Izvestiya,
vol.~16,
2019,
pp.~435--438,
doi:10.33048/semi.2019.16.025.

\bibitem{OstrovskyNafpaktos2023}
A.~Ostrovsky,
\textit{Finitely open or closed functions},
Conference on Topology and its Applications,
Nafpaktos, Greece,
2023,
p.~114.



\bibitem{Hertling1996}
P.~Hertling,
\textit{Topological Complexity of Real Functions},
PhD thesis,
FernUniversität Hagen,
1996,

\end{thebibliography}
\end{document}